\documentclass[12pt,oneside]{article}
\usepackage{amsmath,amsthm,amsfonts,amssymb,MnSymbol,mathrsfs,latexsym}
\usepackage{mathtools}
\usepackage[english]{babel}

\usepackage{blindtext}
\usepackage{geometry}

\newcommand{\const}{\rm const}

\DeclareMathOperator*{\esssup}{ess\,sup}

\theoremstyle{plain}
\newtheorem{theorem}{Theorem}[section]

\newtheorem{proposition}[theorem]{Proposition}
\newtheorem{definition}{Definition}[section]
\newtheorem{remark}{Remark}[section]

\renewenvironment{proof}{{\bf{Proof.}}}{\hfill $\Box$ \\}

\title{\large \textbf{Relations between Lorentz-Zygmund\\ and Grand Lebesgue Spaces and  norms}}

\footnotesize\date{}

\author{\normalsize Maria Rosaria Formica ${}^{1^*}$,   \normalsize Eugeny Ostrovsky
${}^2$ and \normalsize Leonid Sirota ${}^3$}

\begin{document}

\maketitle

\begin{center}
{\footnotesize ${}^{1^*}$ Universit\`{a} degli Studi di Napoli \lq\lq Parthenope\rq\rq, via Generale Parisi 13,\\
Palazzo Pacanowsky, 80132,
Napoli, Italy.} \\

\vspace{2mm}

{\footnotesize e-mail: mara.formica@uniparthenope.it} \\

\vspace{4mm}

{\footnotesize ${}^{2,\, 3}$  Bar-Ilan University, Department of Mathematics and Statistics, \\
52900, Ramat Gan, Israel.} \\

\vspace{2mm}

{\footnotesize e-mail: eugostrovsky@list.ru}\\

\vspace{2mm}

{\footnotesize e-mail: sirota3@bezeqint.net} \\

\end{center}

\vspace{0.1cm}

\begin{abstract}
 We establish imbedding properties between Grand Lebesgue Spaces and (generalized) Lorentz-Zygmund ones.
We extend some known previous results concerning imbedding theorems
between Grand Lebesgue and classical Lebesgue-Riesz spaces and we
show the exactness of the obtained estimates.
\end{abstract}

\vspace{5mm}


\vspace{2mm}


 \vspace{5mm}

\section{Introduction.}

\vspace{4mm}

\subsection {\it Classical Lorentz-Zygmund-Orlicz spaces.}

\vspace{4mm}

 \hspace{3mm} Let $ \ (\Omega = \{\omega\}, \cal{B}, {\bf P})  \ $ be a non-trivial probability space with expectation $ \ {\bf E}. \ $
Define, as usually, for an arbitrary numerical valued random
variable (r.v.) $ \  \xi: \ \Omega \to \mathbb R \ $  its tail
function

$$
T_{\xi}(t) = T[\xi](t) \stackrel{def}{=} {\bf P} ( |\xi| > t), \ t \ge 0.
$$

Introduce the following classical Young-Orlicz function

\begin{equation} \label{Young Orl fun}
N_{p,\alpha}(u) \stackrel{def}{=} |u|^p \ [\ln(e + |u|)]^{\alpha}, \
\ p \ge 1, \ \alpha \in \mathbb R, \ \ u \in \mathbb R,
\end{equation}
and the corresponding Luxemburg (Young-Orlicz) norm for the random
variable (r.v.) $ \ \xi, \ $, as follows
\begin{equation}\label{Young norm}
\begin{split}
||\xi||_{LN_{p,\alpha}}& =   ||\xi||_{p,\alpha}
\stackrel{def}{=}\inf \{\mu > 0: \ \int_{\Omega} |\xi(\omega)/\mu|^p
\ \ln^{\alpha}(e +
|\xi|/\mu) \ {\bf P}(d \omega) \le 1 \ \}\\
& = \inf \{\mu > 0: \  {\bf E} |\xi(\omega)/\mu|^p \ \ln^{\alpha}(e
+ |\xi|/\mu) \le 1 \}.
\end{split}
\end{equation}

The corresponding Orlicz space
$$
L[p,\alpha] \stackrel{def}{=} \{ \ \xi: \ ||\xi||_{p,\alpha}<\infty\
\}
$$
 is named {\it Lorentz-Zygmund space} or, equally, {\it Lorentz-Zygmund-Orlicz} space. \par
 \ These spaces were applied, in particular, in functional analysis (\cite{Bennett Rudnick,Bennet Sharpley,
 Edmunds,Lorentz1,Lorentz2});
 \ theory of extrapolation of operators (\cite{Fiorenza Krbec}), theory of trigonometric series (\cite{Zygmund}),
theory of probability (\cite{Braverman1991}), theory of
approximation (\cite{Ostr Nikol type}), etc.

\vspace{3mm}

 \ Of course, the space $ \  L[p, 0]  \ $ coincides with the ordinary Lebesgue-Riesz one $  \  L[p] \ $ builded on the
 source probability triplet. \par

 \vspace{3mm}

  \hspace{3mm} Let us study the relations between Lorentz-Zygmund norms of a random variable and its tail of distribution behavior.
 Suppose first of all that for some non - negative r.v. $ \eta \ $ there holds

$$
||\eta||_{LN_{p,\alpha}} =   ||\eta||_{p,\alpha} = K \in (0,\infty);
$$
  then

 $$
 {\bf E} N_{p,\alpha}(\eta/K) \le 1.
 $$
  \ It follows from Tchebyshev's - Markov's inequality

 \begin{equation} \label{tail prlim est}
{\bf P} (\eta/K > t) \le 1/N_{p,\alpha}(t), \ t > 0,
 \end{equation}
therefore

 \begin{equation} \label{tail est}
{\bf P} (\eta/K > t) \le \Delta_{p,\alpha}(t), \ t > 0,
 \end{equation}
where by definition

 \begin{equation} \label{Delta p alpha}
\Delta_{p,\alpha}(t) \stackrel{def}{=} \min[1, 1/N_{p,\alpha}(t)], \ t > 0,
 \end{equation}
so that  the function $ \ t \to \Delta_{p,\alpha}(t), \ t \ge 0 \ $,
is actually some  {\it tail function,}  i.e.  it is the tail
function for some non - negative r.v.\par

 \hspace{3mm} On the other words, one has under the same conditions and notations, in particular for non - negative r.v.
 $ \ \eta \ $,
\begin{equation} \label{ K tail estim}
T_{\eta}(t) = {\bf P} (\eta > t) \le \Delta_{p,\alpha}(t/K), \ t > 0.
 \end{equation}

 \vspace{4mm}

 \ Conversely,   let the estimate (\ref{tail est}) be given for the non - negative r.v. $ \eta. \ $ One can assume without loss of generality
 $ \ K = 1. \ $ As long as

$$
{\bf E} \eta^s = s \int_0^{\infty} t^{s-1} \ {\bf P}(\eta > t) \ dt,  \ s > 0,
$$
following  we have in the considered case, denoting $ \ \delta =
p-s, \ \delta \in (0, s/p) \ $,

\begin{equation}\label{Expectation J estim}
{\bf E} |\eta|^s = ||\eta||_s^s \le J(\alpha, p,s),
\end{equation}
where
\begin{equation}  \label{J estim}
 J(\alpha, p,s) \stackrel{def}{=} \  s \int_0^{\infty} \left[ \ t^{s-1} \cdot \min [ \ 1, 1/N_{p,\alpha}(t) \ ] \ \right]  dt,  \ s > 0.
\end{equation}

\vspace{3mm}

\begin{proposition}
For the values $s \in [0,p) \ $ we have
$$
J(\alpha, p,s) \le s \cdot (C_1(\alpha,p) + V(\alpha, \delta)),
$$
where
\begin{equation*}
\begin{split}
\alpha < 1 & \ \ \Rightarrow \  \ V(\alpha, \delta) \asymp
C_2(\alpha)
\delta^{\alpha - 1},\\ \\
\alpha = 1 & \ \ \Rightarrow \ \ V(1, \delta) \asymp C_3 \ (1 + |\ln
\delta|),\\ \\
\alpha > 1 & \ \ \Rightarrow \ \ V(\alpha, \delta) \asymp
C_4(\alpha) \in (0,\infty),
\end{split}
\end{equation*}
where $V, C_1, C_2, C_3, C_4$ are positive constants.
 The corresponding norm estimations in the Lebesgue-Riesz space $L_{s}$ also hold.
\end{proposition}

\vspace{4mm}

 \begin{remark}
 {\rm The estimate \eqref{Expectation J estim} is essentially non-improvable, for instance for the positive r.v.
 $ \ \eta \ $ such that
 $ \ \exists K =\const \in (0,\infty) \ \Rightarrow \ $

 \begin{equation} \label{tail beta est}
{\bf P} (\eta/K > t) \asymp \min[1, 1/N_{p,\alpha}(t)], \ t > 0.
 \end{equation}
 }
 \end{remark}

 \vspace{4mm}

\subsection{\it Grand Lebesgue Spaces (GLS)}

\vspace{4mm}

 \ Let $p\geq 1 $. The usual Lebesgue - Riesz norm $ \ ||f||_p = ||f||_{L_p(\Omega)} =||f||_{L_p({ \Omega, \bf P})} \ $
 of a measurable function
 $ \ f: \Omega \to \mathbb R \ $ is defined as
$$
||f||_p =   \left[ \ {\bf E} |f|^p \ \right]^{1/p} =
||f||_{L_p(\Omega,{\bf P})} \stackrel{def}{=} \left\{  \int_{\Omega}
\ |f(\omega)|^p \ {\bf P}(d \omega) \  \right\}^{1/p} = \left[ \
I(|f|^p) \ \right]^{1/p}, \ \ 1\leq p <\infty,
$$
and of course
$$
||f||_{\infty} = \esssup_{\omega \in \Omega}|f(\omega)|, \ \ \
p=\infty.
$$

\vspace{4mm}

 \ We recall here for reader convenience some known definitions and  facts  from the theory of Grand Lebesgue Spaces (GLS).
    \ Let $ 1 \le  a < b \le \infty$ and $ \ \psi = \psi(p), \ p \in (a,b)$, be a positive measurable numerical valued
 function, not necessary finite in $a,b$, such that $ \displaystyle \inf_{p \in (a,b)} \psi(p) > 0. \ $

We denote with $ G \Psi(a,b)$ the set of all such functions
$\psi(p), \ p \in (a,b)$, for some  $ 1 \le a < b \le \infty \ $ and
put
$$
G\Psi := \displaystyle\bigcup_{1 \le a < b \le \infty} G \Psi(a,b).
$$
For instance, if $m = {\const} > 0$,
$$
\psi_m(p) := p^{1/m},  \ \ p \in [1,\infty)
$$
or, for $ 1 \le a < b < \infty, \ \alpha,\beta = \const \ge 0$,
$$
   \psi^{(a,b; \alpha,\beta)}(p) := (p-a)^{-\alpha} \ (b-p)^{-\beta}, \  \ p \in
   (a,b),
$$
the above functions belong to $ G \Psi(a,b)$.

The (Banach) Grand Lebesgue Space$ \ G \psi  = G\psi(a,b)$ consists
of all the real (or complex) numerical valued measurable functions
$f: \Omega \to \mathbb R$ having finite norm
\begin{equation} \label{norm psi}
    || \ f \ || = ||f||_{G\psi} \stackrel{def}{=} \sup_{p \in (a,b)} \left[ \frac{||f||_p}{\psi(p)} \right].
 \end{equation}
We write $ \ G\psi \ $ when $ \ a = 1\ $ and $ \ b= \infty. \ $ The
function $ \  \psi = \psi(p) \  $ is named the {\it generating
function } for the space $G \psi$ . \par
If for instance
$$
  \psi(p) = \psi^{(r)}(p) = 1, \ \  p = r;  \  \ \psi^{(r)}(p) = +\infty,   \ \ p \ne r,
$$
 where $ \ r = {\const} \in [1,\infty),  \ C/\infty := 0, \ C \in \mathbb R, \ $ (an extremal case), then the correspondent
 $ \  G\psi^{(r)}(p)  \  $ space coincides  with the classical Lebesgue - Riesz space $ \ L_r = L_r(\Omega). \ $ \par

\vspace{4mm}

These spaces are investigated in many works (e.g.
\cite{KozOs,Lifl,Ostrovsky2,Ostrovsky4,Ostrovsky5}). For example
they play an important role in the theory of Partial Differential
Equations (PDEs) (see,
e.g.,\cite{Greco-Iwaniec-Sbordone-1997,Fiorenza-Formica-Gogatishvili-DEA2018,fioformicarakodie2017,
Ahmed-Fiorenza-Formica-Gogatishvili-Rakotoson}), in interpolation
theory (see, e.g.,
\cite{fioforgogakoparakoNA,fiokarazanalanwen2004}), in the theory of
Probability (\cite{Ermakov etc.
1986,Ostrovsky3,Ostrovsky5,ForKozOstr_Lithuanian}), in Statistics
\cite[chapter 5]{Ostrovsky1}, in theory of random fields
\cite{KozOs,Ostrovsky4}.

These spaces are rearrangement invariant (r.i.) Banach function
spaces; the fundamental function has been studied in
\cite{Ostrovsky4}. They not coincide, in the general case, with the
classical Banach rearrangement functional spaces: Orlicz, Lorentz,
Marcinkiewicz, etc., (see \cite{Lifl}, \cite{Ostrovsky2}). The
belonging of a function $ f: \Omega \to \mathbb{R}$ to some $ G\psi$
space is closely related to its tail function behavior
$$
 T_f(t) \stackrel{def}{=} {\bf P}(|f| \ge t), \ \ t \ge 0,
 $$
as $ \ t \to 0+ \ $ as well as when $ \ t \to \infty $ (see
\cite{KozOs,KozOsSir2017,KozOsSir2019}).

\vspace{4mm}

\subsection{\it Generalized Lorentz-Zygmund Spaces}

\vspace{4mm}

\hspace{3mm} For each measurable function $ f: \Omega \to \mathbb R$
we denote, as usually, by $ \ d_f \ $  and $ \ f^* \ $ the {\it anti
- distribution function} and the non-increasing rearrangement of
$f$, respectively:
$$
 d_f(\sigma) := {\bf P} \{\omega \in \Omega \, : \, |f(\omega)| \ge \sigma \}, \hspace{5mm}
 f^*(t) := \inf \{\sigma > 0 \, : \, d_f(\sigma) < t \ \}.
$$
 \ Evidently,

$$
\int_0^{\infty} [f^*(t)|^p \ dt = \int_{\Omega} |f(\omega)|^p \ {\bf
P}(d \omega) ={\bf E}|f|^p, \ \ p > 1.
$$

\vspace{3mm}

 \ As a slight consequence: for an arbitrary generating function $ \ \psi(\cdot) \in\ G\Psi \ $
 and for an arbitrary random variable
 $ \ f: \ \Omega \to \mathbb R \ $ there holds

\begin{equation} \label{coincidence}
||f^*||_{G\psi (\mathbb R_+)}  = ||f||_{G\psi(\Omega)}.
\end{equation}

\vspace{4mm}

Let $S = S(t), \, t \in (0,1)$, be a numerical valued non-negative
integrable function. Define the following {\it generalized
Lorentz-Zygmund norm} \ $ ||f||_{[S]} \ $ for an arbitrary
measurable function $ \ f: (0,1) \to \mathbb R \ $
\begin{equation} \label{S norm}
||f||_{[S]} \stackrel{def}{=} \int_0^1 f^*(t) \ S(t) \ dt.
\end{equation}
It  follows from H\"older's inequality
\begin{equation} \label{S Holder}
||f||_{[S]} \le ||f||_p \cdot ||S||_{p'}, \ \ p \in (1, \infty), \ \
p' := p/(p-1).
\end{equation}

\vspace{3mm}

 \ Assume now that the function $ f(\cdot)$ belongs to some $ G\psi(a,b) $ space and $S(\cdot)$ belongs to some $G\nu(c,d)$
 space, with $\psi(\cdot)$ and $\nu(\cdot)$ generating functions. Then it follows from \eqref{S Holder}
$$
||f||_{[S]} \le \{ \ ||f||_{G\psi(a,b)} \, \cdot \,
||S||_{G\nu(c,d)} \ \} \cdot \psi_{a,b}(p) \cdot \nu_{c,d}(p').
$$

Hence we have the following

\vspace{3mm}

\begin{proposition}

\begin{equation} \label{Altogeth}
||f||_{[S]} \le  \zeta \cdot ||f||_{G\psi(a,b)} \cdot
||S||_{G\nu(c,d)},
\end{equation}

where

\begin{equation} \label{nuu}
\zeta = \zeta(a,b; c,d) := \inf_{ p \in (a,b); \ p' \in (c,d)} \left\{ \ \psi_{a,b}(p) \cdot \nu_{c,d}(p')  \ \right\}.
\end{equation}
The \lq\lq constant\rq\rq $ \zeta = \zeta(a,b; c,d)$ in
\eqref{Altogeth}, \eqref{nuu} is the best possible (see e.g.
\cite{Lifl,Ostrovsky5}). As ordinarily, $ \ \inf _{p \in \emptyset}
(\cdot) = + \infty. \ $  \par

\end{proposition}

\vspace{4mm}

\subsection{\it Grand Zygmund Spaces}

\vspace{4mm}

Let us consider a certain sub-domain $ \ Q \ $ on the quarter plane
\begin{equation} \label{quarter}
Q  \subset\{ (1,\infty) \otimes (0,\infty) \}.
\end{equation}
Let also $ \ \rho(p,\gamma), \ (p,\gamma) \in Q \ $, be some
positive continuous numerical valued function.

\vspace{4mm}

\begin{definition}

\vspace{3mm}

{\rm \ The so-called Grand Zygmund Spaces (GZS) $ \ GZ[Q](\rho) \ $
consists of all the measurable functions (random variables)
 $ \ \{\eta\} \ $ having a finite norm}

\begin{equation} \label{GZS norm}
V = V(\eta) := ||\eta||_{GZ[Q](\rho)}  \stackrel{def}{=}
\sup_{(p,\gamma) \in Q}  \left[ \
\frac{||\eta||_{p,\gamma}}{\rho(p,\gamma)} \ \right].
\end{equation}
\end{definition}

 \vspace{4mm}

 \hspace{3mm} {\bf Our aim in this  short report is to establish the interrelation properties between tail behavior of a
 random variable and its belonging to appropriate (generalized) Grand Zygmund Space}.

\vspace{3mm}

\begin{remark}
{\rm Let $ \ \nu(p) = \rho(p,0); \ $ then  the space $ \ GZ[Q](\nu)
\ $ coincides with the Grand Lebesgue Space $ \ G\nu. \ $
 This case was investigated in many works (\cite{Buld Koz AMS,Ermakov etc.
 1986,KozOs,KozOsSir2017}).
}
\end{remark}

 \vspace{3mm}

 \ Notice that the embedding Theorems between generalized weighted Lorentz spaces are investigated in the recent article
 \cite{Gogatishvili at all}. \par

\vspace{3mm}

\begin{proposition}
Suppose that for some random variable $ \ \eta \ $
\begin{equation} \label{GZS tail estim}
V = V(\eta) := ||\eta||_{GZ[Q](\rho)}  \in (0,\infty).
\end{equation}
Then one has the following tail estimate
\begin{equation} \label{tail behav est}
T_\eta(t) \le \inf_{(p,\gamma) \in Q} \Delta_{p,\gamma}[t/(V\cdot
\rho(p,\gamma))],
\end{equation}
where $\Delta_{p,\gamma}$ is defined as in \eqref{Delta p alpha}.
\end{proposition}

 \vspace{3mm}

\begin{proof}
 The proof is very simple. Let $ \ (p,\gamma) \ $ be an arbitrary point from the set $ \ Q; \ $
one can apply the inequality \eqref{ K tail estim}
\begin{equation} \label{K gamma estim}
T_{\eta}(t) = {\bf P} (\eta > t) \le \Delta_{p,\gamma}(t/(V \cdot \rho(p,\gamma))), \ t > 0.
 \end{equation}

 \ It remains to take the minimum  in (\ref{K gamma estim})  over both the parameters $ \ (p,\gamma) \ $ from the set $ \ Q. \ $ \par
\end{proof}

 \vspace{1mm}

\begin{remark}
{\rm The last estimate (\ref{K gamma estim}) is not improvable yet
in the case when $ \ \gamma = 0, \ $ as long as this case
 coincides with the well studied case of Grand Lebesgue Spaces. }
 \end{remark}

 \vspace{4mm}

\section{Main result: estimates}

\vspace{4mm}

 \begin{center}

 {\it Auxiliary facts.}

 \end{center}

\vspace{4mm}

 \  {\sc 1. Dilation of Young - Orlicz function. } \par

   \vspace{3mm}
   \ Define for an arbitrary Young - Orlicz function $ \ \Phi = \Phi(u), \ u \in \mathbb R$, its dilation transform

 \begin{equation} \label{dilat transf}
 T_C \Phi(u) := \Phi(Cu), \ \ C = \const > 0.
 \end{equation}

  One has, for an arbitrary r.v. $ \ \xi \ $,

\begin{equation} \label{dilation  equality}
||\xi||_{L(T_C \Phi)}  = C ||\xi||_{L\Phi}.
\end{equation}

\vspace{4mm}

  \ {\sc 2. Monotonicity of Orlicz's norm.}

\vspace{3mm}

 Let $ \  \Phi_1(\cdot), \ \Phi_2(\cdot)   \ $ be two Young - Orlicz functions such that
$$
 \Phi_1(u) \le  \Phi_2(u), \ \ \ \forall u \in \mathbb R .
 $$
\ Then, for any r.v. $ \xi$,
\begin{equation}  \label{monotonicity}
||\xi||_{L\Phi_1}  \ge ||\xi||_{L\Phi_2}
\end{equation}
\ holds.

 \vspace{3mm}

 \ For instance, if $ \ 1 \le r \le p \le \infty, \ $ then $ \ ||\xi||_r \le ||\xi||_p$ \, (Jensen's inequality). \\
 \ Recall that our measure $ \ {\bf P} \ $ is probabilistic. \par

%
%

\vspace{4mm}

\subsection{Upper estimate}

\vspace{4mm}

 \hspace{3mm}  Define the function

\begin{equation} \label{glav fun}
\Theta(r,p,\gamma) \stackrel{def}{=} \left\{ \Gamma [ (\gamma p)/(p - r) + 1 ] \right\}^{(p-r)/pr}, \ \gamma \ge 0,  \ 1 \le r < p < \infty.
\end{equation}
From some detailed calculations in the articles \cite{Cap Fior
Krbec,Fiorenza Krbec} it follows

\begin{proposition}

\begin{equation} \label{main estimate}
||g||_{r,\gamma} \le \Theta(r,p,\gamma) \cdot ||g||_p.
\end{equation}
In particular,
\begin{equation} \label{part estim}
||g||_{1,\gamma} \le \Theta(1,p,\gamma) \cdot ||g||_p, \ \gamma \ge 0, \ p > 1.
\end{equation}
\end{proposition}

\vspace{4mm}

\begin{remark}

{\rm  Our estimates \eqref{main estimate} and, consequently,
\eqref{part estim} are exact. Namely, denote

$$
U[\gamma,r,p](f) := \int_0^1 f^*(t) \cdot |\ln t|^{\gamma} \ dt,
$$

and
\begin{equation} \label{exactness}
K = K[\gamma, r,p] \stackrel{def}{=} \sup_{0 \ne f \in L(p)} \left[ \ \frac{|U[\gamma,r,p](f)|}{||f||_p} \ \right],
\end{equation}
where as above $ \ \gamma \ge 0, \ 1 \le r < p < \infty. \ $  We assert

\begin{equation} \label{unimprovability}
K[\gamma, r,p] =  \Theta(r,p,\gamma).
\end{equation}
 \  Truly, the {\it upper} estimate is contained in (\ref{main estimate});  as well as the {\it lower one} is attained on
 the test function
$$
f_0(t) := |\ln t|^{\gamma p/(p-1)}, \ t \in (0,1),
$$
or equivalently, $ \ g_0(t) = C \ f_0(t), \ C = \const \in
(0,\infty). \ $ We applied the well - known  conditions for the {\it
lower} estimate in H\"older's inequality. }
\end{remark}

\vspace{4mm}

\subsection{Lower estimate}

\vspace{4mm}

 \hspace{3mm} We use the following inequality (in our notations)

\begin{equation} \label{aux ps}
||g||_{p, - \beta} \ge C(\beta) \ \left[ \ \frac{s}{p-s}  \
\right]^{-\beta/s} \ ||g||_s,
\end{equation}
where  $ \ \beta > 0, \  C(\beta) \in (0,\infty), \ $

$$
 p \in [1,\infty), \ s \in [1,p);
$$
the case $ \ \beta = 0 \ $ is trivial. See  \cite{Fiorenza
Krbec,Gogatishvili at all,Ostr Nikol type}.

 \ Introduce the following function

\begin{equation} \label{aux inver fun}
\kappa[g](s) \stackrel{def}{=} \inf_{p \ge s} \ \left\{ \ \left[\frac{s}{p- s} \right] \cdot \ ||g||_{p, - \beta} \right\}.
\end{equation}

\vspace{4mm}

\begin{proposition}
Suppose that the function $ \ g(\cdot) \ $ is such that for some $ \
(a,b), \ 1 \le a < b \le \infty $, it is $ \ \kappa[g](s) < \infty.
\ $ Then

\begin{equation} \label{inver est}
||g||_s \le C(\beta) \ \kappa[g](s), \ \  s \in (a,b).
\end{equation}

%
On the other words,

\begin{equation} \label{GPsi inver}
||g||_{G\upsilon} \le C(\beta),  \hspace{3mm} \hbox{\rm where}
\hspace{3mm} \upsilon(s) := \kappa[g](s).
\end{equation}

\end{proposition}

 \vspace{5mm}

\section{Concluding remarks.}

\vspace{4mm}

 \hspace{3mm} It is interesting in our opinion to investigate more general Orlicz  spaces having the correspondent
 Young - Orlicz function of the form

$$
N[p,S](u) \stackrel{def}{=} |u|^p \ S(|u|),
$$
 where $ \ p \ge 1, \ S = S(u) \ $ is non - negative even continuous {\it slowly varying} at infinity function,
 so that the function $ \ u \to N[p,S](u) \ $ is {\it regular varying}. \par

\vspace{6mm}

\emph{\textbf{Acknowledgements}.} {\footnotesize 
M.R. Formica is member of Gruppo
Nazionale per l'Analisi Matematica, la Probabilit\`{a} e le loro
Applicazioni (GNAMPA) of the Istituto Nazionale di Alta Matematica
(INdAM) and member of the UMI group \lq\lq Teoria
dell'Approssimazione e Applicazioni (T.A.A.)\rq\rq }

\end{document}